\documentclass[12pt,leqno]{amsart}
\usepackage{amssymb}
\usepackage{amsmath,amssymb,color}
\numberwithin{equation}{section}

\newcommand{\weight}{e^{2s\varphi}}
\newcommand{\tilweight}{e^{2s\widetilde{\va}}}
\newcommand{\ep}{\varepsilon}
\newcommand{\la}{\lambda}
\newcommand{\va}{\varphi}
\newcommand{\ppp}{\partial}

\newcommand{\www}{\widetilde}

\newcommand{\R}{\mathbb{R}}

\newcommand{\ooo}{\overline}
\newcommand{\OOO}{\Omega}

\newcommand{\wdel}{\www{\delta}}
\newcommand{\weps}{\www{\varepsilon}}

\allowdisplaybreaks
\title
[]
{
Stability for inverse source problems by Carleman estimates
}

\author{
$^1$ X.~Huang,
$^2$ O.~Yu.~Imanuvilov and \, $^{3,4,5}$ M.~Yamamoto }

\thanks{
$^1$ Graduate School of Mathematical Sciences, The University of Tokyo, 
Komaba, Meguro, Tokyo 153-8914, Japan 
e-mail: {\tt huangxc@ms.u-tokyo.ac.jp}\\
$^2$ Department of Mathematics, Colorado State
University, 101 Weber Building, Fort Collins, CO 80523-1874, U.S.A.
e-mail: {\tt oleg@math.colostate.edu}\\
$^3$ Graduate School of Mathematical Sciences, The University
of Tokyo, Komaba, Meguro, Tokyo 153-8914, Japan \\
$^4$ Honorary Member of Academy of Romanian Scientists, 
Splaiul Independentei Street, no 54,
050094 Bucharest Romania \\
$^5$ Peoples' Friendship University of Russia 
(RUDN University) 6 Miklukho-Maklaya St, Moscow, 117198, Russian Federation
e-mail: {\tt myama@ms.u-tokyo.ac.jp}
}

\date{}
\begin{document}
\maketitle

\baselineskip 18pt

\begin{abstract}
In this article, we provide a modified argument for proving conditional 
stability for inverse problems of determining spatially varying 
functions in evolution equations by Carleman estimates.
Our method needs not any cut-off procedures and can simplify the 
existing proofs.
We establish the conditional stability for inverse source problems 
for a hyperbolic equation and a parabolic equation, and our method 
is widely applicable to various evolution equations.
%\end{abstract}
\\
{\bf Key words.}  
inverse source problem, Carleman estimates, stability
\\
{\bf AMS subject classifications.}
35R30, 35R25
\end{abstract}

\section{Introduction and main results}
For evolution equations, we consider inverse source problems of determining 
spatially varying functions in non-homogeneous terms of the equations.  

More precisely, let $\OOO \subset \R^n$ be a bounded open domain
with smooth boundary $\ppp\OOO$, and let $x=(x_1,..., x_n)\in \R^n$ and 
$t$ denote the spatial and the time variables respectively.  We set
$$
\ppp_{x_j} = \frac{\ppp}{\ppp x_j}, \quad
\ppp_{x_ix_j}^2 = \frac{\ppp^2}{\ppp x_j\ppp x_i}, \quad 1\le j \le n, \quad 
\ppp_t = \frac{\ppp}{\ppp t}, \quad \nabla = (\ppp_{x_1}, ..., 
\ppp_{x_n}), 
$$
$$
\nabla_{x,t} = (\nabla, \ppp_t), \quad \Delta = \sum_{j=1}^n \ppp_{x_j}^2.
$$
By $\nu = \nu(x)$ we denote the unit outward normal vector to 
$\ppp\OOO$ at $x$, and set $\ppp_{\nu}u = \nabla u \cdot \nu$. 

Let $\mathcal{L}$ be a suitable partial differential operator in $(x,t)$
and $I$ be an open time interval.
We consider
$$
\mathcal{L}u = R(x,t)f(x), \quad x\in \OOO, \, t\in I.     \eqno{(1.1)}
$$

Our inverse problem is formulated as follows:
\\
{\it For given $t_0 \in I$, function $R(x,t)$ and subboundary 
$\Gamma \subset \ppp\OOO$, determine $f(x)$ in (1.1) by 
$u\vert_{\Gamma\times I}$, $\nabla u\vert_{\Gamma\times I}$ and
$u\vert_{t=t_0}$.}
\\

The choices of the operator $\mathcal{L}$ in (1.1) are quite general, and 
typical cases are 
$$
\mathcal{L}u = \ppp_t^ku - \Delta u - \sum_{j=1}^n b_j(x)\ppp_{x_j}u - c(x)u,
\quad \mbox{$k=1$ or $k=2$}                      \eqno{(1.2)}
$$ 
with $b_j, c \in L^{\infty}(\OOO)$, $j=1, ..., n$.
We can similarly consider more general elliptic operators but here we omit.

Our formulation for the inverse problem 
requires only a single measurement of data of 
solution to an initial boundary value problem for (1.1).
For our inverse problem, Bukhgeim and Klibanov \cite{BK} created a fundamental 
methodology which is based on Carleman estimates, and established 
the uniqueness for inverse problems.  See also Klibanov \cite{Kl2}, 
\cite{Kl2013}.

A Carleman estimate is an $L^2$-weighted estimate for solutions to 
system (1.1), and is stated as follows: 
by choosing a weight function $\va = \va(x,t)$,
there exist constants $C>0$ and $s_0>0$ such that 
$$
\int_{\OOO\times I} s^3\vert u\vert^2 \weight dxdt
\le C\int_{\OOO\times I} \vert \mathcal{L}u\vert^2 \weight dxdt
+ C\int_{\ppp(\OOO\times I)} \vert \nabla_{x,t}u\vert^2 \weight d\Sigma
                                           \eqno{(1.3)}
$$
for all $s\ge s_0$.  We note that the constant $C>0$ should be independent of 
$s \ge s_0$.  The choices of the weight function $\va(x,t)$ are essential
for the applications, and in this paper we use two types of weight functions:
$$
\va(x,t) = e^{\lambda(d(x) - \beta (t-t_0)^2)}   \eqno{(1.4)}
$$
and
$$
\va(x,t) = \exp\left( 
\frac{ e^{\lambda d(x)} 
- e^{2\lambda \Vert d\Vert_{C(\ooo{\OOO})}} }{t(T-t)}\right),   \eqno{(1.5)}
$$
where $\lambda>0$ is a large constant, and $d$ is a suitable 
function.  Carleman esimates with the weight function (1.5) hold for 
parabolic and Schr\"odinger equations (Imanuvilov \cite{Ima1995}, 
Imanuvilov and Yamamoto \cite{IY1998}, Baudouin and Puel \cite{BP}),
but not for hyperbolic types of equations, while the ones with (1.4) more
comprehensively hold.

Since \cite{BK}, we have had many works on inverse problems on the 
basis of Carleman estimates.
Among them, Imanuvilov and Yamamoto \cite{IY1998},
\cite{IY1}, \cite{IY2} are early works establishing 
the best possible Lipschitz stability over the whole domain $\OOO$.  

As monographs, we can refer to 
Beilina and Klibanov \cite{BeiKl}, Bellassoued and Yamamoto \cite{BY2017},
Fu, L\"u and Zhang \cite{FLZ}, Klibanov and Timonov \cite{KlT}.
Moreover we list some of related articles on inverse problems by 
Carleman estimates.  Since the researches have been 
developing widely, it is not easy to compose any comprehensive 
lists, and one can also consult the references therein.

{\bf Hyperbolic equations}.\\
Beilina, Cristofol, Li and Yamamoto \cite{BCLY},
Imanuvilov and Yamamoto \cite{IY3}.  

{\bf Parabolic equations}.\\
Imanuvilov and Yamamoto \cite{IY1998}, Yamamoto and Zou \cite{YZ}
by Carleman estimates of type (1.5), and Yamamoto \cite{Ya} as survey. 

{\bf Schr\"odinger equations}.\\
Baudouin and Mercado \cite{BM}, Baudouin and Puel \cite{BP},
Mercado, Osses and Rosier \cite{MOR} by Carleman estimates of type (1.5), and
Yuan and Yamamoto \cite{YY2} by Carleman estimates of type (1.4).

{\bf First-order equations (transport equations)}.\\
Cannarsa, Floridia and Yamamoto \cite{CFY},
Cannarsa, Floridia, G\"olgeleyen and Yamamoto \cite{CFGY},
G\"olgeleyen and Yamamoto \cite{GY2016}. 

For plate equations and integro-differential equations
related to the viscoelasticity, see for example, 
Yuan and Yamamoto \cite{YY1}, Cavaterra, Lorenzi and Yamamoto \cite{CLY}.
\\

In the existing works, whenever one applied Carleman estimates of type 
(1.4), one needed to 
introduce cut-off functions $\chi(t)$ or $\chi(x,t)$ in order that 
$\chi u$ vanishes on the boundary of the domains in $x$ and $t$ where 
we do not know data of $u$.  On the other hand, in applying Carleman estimates
of type (1.5), we need not any cut-off.

The cut-off procedure makes the arguments for the inverse problems 
more complicated, because we have to apply Carleman estimates not directly to 
solution to (1.1), but to the functions multiplied by $\chi$, 
and the structure of the original equations may be 
changed inconveniently.

In this article, we propose an argument without the cut-off procedure for
proving the stability for the inverse problems on the basis of 
Carleman estimates of type (1.4).
The key is that the weight function already takes smaller values on 
the boundary of a domain in $(x,t)$ where data are not given, so that 
the weight function can well control such unknown data 
for proving the stability in the inverse problems, and so the cut-off is 
not necessary.

Here we investigate an inverse source problem  only for  second order 
differential operators $\mathcal L$  of  hyperbolic or  parabolic types  but
our argument can work similarly to other evolution equation which admit 
suitable apriori estimates of Carleman type.
\\

First we consider an inverse source problem for a hyperbolic equation.
$$
\left\{ \begin{array}{rl}
& \ppp_t^2u - \Delta u - \sum_{j=1}^n b_j(x)\ppp_{x_j}u - c(x)u
= R(x,t)f(x), \quad x \in\OOO, \, 0<t<T, \\
& u(\cdot,0) = \ppp_tu(\cdot,0) = 0 \quad \mbox{in $\OOO$},\\
& u\vert_{\ppp\OOO\times (0,T)} = 0.
\end{array}\right.                         \eqno{(1.6)}
$$
Here we assume $b_j, c \in L^{\infty}(\OOO)$, $j=1,..., n$.

For arbitrarily fixed $x_0 \not\in \ooo{\OOO}$, we set 
$$
\Gamma := \{ x\in \ppp\OOO;\, (x-x_0) \cdot \nu(x) \ge 0\}.
                                                       \eqno{(1.7)}
$$
We can prove the following.
\\
{\bf Theorem 1 (global Lipschitz stability for an inverse 
source problem for a hyperbolic equation)}.\\
{\it We assume that there exists a constant $r_0 > 0$ such that 
$$
R \in H^1(0,T;L^{\infty}(\OOO)), \quad \vert R(x,0)\vert \ge r_0, 
\quad x \in \ooo{\OOO}                \eqno{(1.8)}
$$
and
$$
T > \left( \max_{x\in \ooo{\OOO}} \vert x-x_0\vert^2
- \min_{x\in \ooo{\OOO}} \vert x-x_0\vert^2\right)^{\frac{1}{2}}.
                                                       \eqno{(1.9)}
$$
Then there exists a constant $C>0$ such that 
$$
\Vert f\Vert_{L^2(\OOO)} \le C\Vert \ppp_t\ppp_{\nu}u\Vert
_{L^2(\Gamma \times (0,T))}
$$
for each $u$ satisfying (1.6) and  the regularity condition
$$
\ppp_tu \in C([0,T];H^2(\OOO) \cap H^1_0(\OOO)) \cap 
C^2([0,T];L^2(\OOO)).                      \eqno{(1.10)}
$$}
\\

We can relax the regularity assumption (1.10) on function $u$, but we omit details for 
simplicity.
This type of stability over $\OOO$ was 
proved by e.g., Imanuvilov and Yamamoto \cite{IY1}
with a cut-off argument in $t$.  The reverse inequality
$$
\Vert \ppp_t\ppp_{\nu}u\Vert_{L^2(\ppp\OOO \times (0,T))}
\le C\Vert f\Vert_{L^2(\OOO)}
$$
can be proved for any $T>0$ by the multiplier method (e.g.,
\cite{BY2017} (Chapter 3), Komornik \cite{Ko}).
\\

Second we consider an inverse source problem for a parabolic 
equation:
$$
\ppp_tu - \Delta u - \sum_{j=1}^n b_j(x)\ppp_{x_j}u - c(x)u
= R(x,t)f(x), \quad x \in\OOO, \, 0<t<T.                       \eqno{(1.11)}
$$
Let $\Gamma \subset \ppp\OOO$ be an arbitrarily fixed non-empty relatively
open subset.  We arbitrarily choose a subdomain $\OOO_0 \subset 
\OOO$ such that $\ooo{\OOO_0} \subset \OOO \cup \Gamma$,
$\ppp\OOO_0 \cap \ppp\OOO$ is a non-empty relatively open subset 
of $\ppp\OOO$ and $\ooo{\ppp\OOO_0 \cap \ppp\OOO} 
\subset \Gamma$.  Let $0 < t_0 < T$ and let 
$I = (t_0 - \delta, t_0 + \delta)$ such that $I \subset (0,T)$.  

Then we have
\\
{\bf Theorem 2 (local H\"older stability for an inverse 
source problem for a parabolic equation)}.\\
{\it We assume that the function $R$ satisfies 
$$
R \in H^1(0,T;L^{\infty}(\OOO)), \quad \vert R(x,t_0)\vert \ge r_0, 
\quad x \in \ooo{\OOO}                \eqno{(1.12)}
$$
for some constant $r_0>0$. Moreover
a pair $(u,f)
\in (H^2(0,T;H^1(\OOO))\cap H^1(0,T;H^2(\OOO))) 
\times L^2(\Omega)$ solve equation (1.11) and 
the function $u$ satisfies an a priori bound:
$$
\Vert u\Vert_{H^2(0,T;H^1(\OOO))} + \Vert u\Vert_{H^1(0,T;H^2(\OOO))}
\le M                               \eqno{(1.13)}
$$
with some constant $M>0$.  
Then there exist constants $C>0$ and $\theta \in (0,1)$ depending 
on $M$, $\Gamma$, $t_0$ such that 
$$
\Vert f\Vert_{L^2(\OOO_0)} 
\le C(\Vert \nabla_{x,t}\ppp_tu\Vert_{L^2(\Gamma\times (0,T))}
+ \Vert \ppp_tu\Vert_{L^2(\Gamma\times (0,T))}
+ \Vert u(\cdot,t_0)\Vert_{H^2(\OOO)})^{\theta}.
$$
}

We note that we have no boundary data on whole $\ppp\OOO\times (0,T)$, but
only $\Gamma \times (0,T)$.
With the whole boundary condition on $\ppp\OOO\times (0,T)$, we can prove
the Lipschitz stability over $\OOO$ by Carleman estimate with type 
(1.5) (Imanuvilov and Yamamoto \cite{IY1998}).
Moreover unlike (1.9) in Theorem 1, we need not any conditions on 
the observation time $T$.

The article is composed of five sections.  In Sections 2, we show the key
Carleman estimates for (1.6) and (1.11).  Sections 3 and 4 are devoted to the 
proofs of Theorems 1 and 2 respectively. Section 5 gives concluding remarks.
\section{Two key Carleman estimates}

We set
$$
Q_{(-T,T)} = \OOO \times (-T,T), \quad Q = \OOO \times (0,T).
$$
We first consider the following hyperbolic equation:
$$
\left\{ \begin{array}{rl}
& \ppp_t^2 v(x,t) - \Delta v(x,t) - \sum_{j=1}^n b_j(x)\ppp_{x_j}v
- c(x)v = F(x,t), \quad (x,t)\in Q_{(-T,T)},\\
& v(x,t) = 0, \quad (x,t)\in \partial\Omega\times (-T,T).
\end{array}\right.
                                               \eqno{(2.1)}
$$

For arbitrarily fixed $x_0 \not\in \ooo{\OOO}$, $\lambda>0$, $0 \le t_0 < T$,
and $0 < \beta < 1$, we set
$$
\va(x,t) = e^{\lambda \psi(x,t)}, \ \psi(x,t) = |x-x_0|^2 - \beta (t-t_0)^2, 
\quad (x,t)\in Q_{(-T,T)}.                     \eqno{(2.2)}
$$

Henceforth $C>0$ denotes generic constants which are independent of 
parameter $s>0$.

{\bf Lemma 1 (Carleman estimate for hyperbolic equation)}.\\
{\it Let $\lambda>0$ be sufficiently large.
Then there exist constants $s_0 > 0$ and $C > 0$ such that 
\begin{align*}
&\int_{Q_{(-T,T)}} \left(s|\nabla_{x,t} v|^2 + s^3|v|^2 \right)e^{2s\va} dxdt 
\le \ C\int_{Q_{(-T,T)}} |F|^2 e^{2s\va} dxdt 
+ C\int_{\Gamma \times (-T,T)} s |\ppp_\nu v|^2 e^{2s\va} d\Sigma\\
+ &C\int_{\Omega} (s|\nabla_{x,t} v(x,T)|^2 + s^3|v(x,T)|^2)e^{2s\va(x,T)} dx\\
+ &C\int_{\Omega} (s|\nabla_{x,t} v(x,-T)|^2 + s^3|v(x,-T)|^2)e^{2s\va(x,0)} dx
\end{align*}
for all $s > s_0$ and $v\in H^2(-T,T; L^2(\Omega))\cap 
L^2(-T,T; H^2(\OOO) \cap H_0^1(\Omega))$ satisfying (2.1). }
\\

Lemma 1 is a classical Carleman estimate and we can prove similarly for example
to Theorem 4.2 in \cite{BY2017} by keeping the values of $u$ at $t=-T, T$.
\\

Second we show a Carleman estimate for a parabolic equation.
We introduce the weight function. 
First we construct some domain $\OOO_1$.
For $\Gamma \subset \ppp\OOO$, we choose a bounded domain $\OOO_1$ 
with smooth boundary such that
$$
\OOO \subsetneqq \OOO_1, \quad \ooo{\Gamma} = \ooo{\ppp\OOO\cap\OOO_1}, 
\quad \ppp\OOO\setminus\Gamma \subset \ppp\OOO_1.
                                                   \eqno{(2.3)}
$$
In particular, $\OOO_1\setminus\ooo{\OOO}$ contains some non-empty 
open subset. 
We note that $\OOO_1$ can be constructed as the interior of a union  
of $\ooo{\OOO}$ and the closure of a non-empty domain 
$\widehat{\OOO}$ satisfying $\widehat{\OOO} \subset \overline{\R^3 
\setminus \OOO}$ and $\ppp\widehat{\OOO} \cap \ppp\OOO = \Gamma$.

We choose a domain $\omega$ such that 
$\ooo{\omega} \subset \OOO_1 \setminus \ooo{\OOO}$.
Then, by \cite{Ima1995}, we can find $d\in C^2(\ooo{\OOO_1})$ such that
$$
d>0 \quad \mbox{in $\OOO_1$}, \quad |\nabla d| > 0 \quad \mbox{on
$\ooo{\OOO_1\setminus \omega}$}, \quad d=0 \quad \mbox{on $\ppp\OOO_1$}. 
                                     \eqno{(2.4)}
$$
In particular, 
$$
d > 0 \quad \mbox{on $\ooo{\OOO_0}$}, \quad 
d=0 \quad \mbox{on $\ppp\OOO \setminus \Gamma$}.    \eqno{(2.5)}
$$
We recall that we choose a domain $\OOO_0 \subset \OOO$ satisfying 
$\ooo{\ppp{\OOO_0} \cap \ppp\OOO} \subset\Gamma$ and 
$\ooo{\OOO_0} \subset \OOO \cup \Gamma$.

Then for arbitrarily fixed $t_0\in (0,T)$ and $\delta>0$ such that 
$0 \le t_0 - \delta < t_0 + \delta \le T$, we set
$$
I = (t_0-\delta, t_0+\delta), \quad Q_I = \OOO \times I.
$$
We define
$$
\www{\psi}(x,t) = d(x) - \beta (t-t_0)^2,\quad 
\www{\varphi}(x,t) = e^{\la\www{\psi}(x,t)}, \quad (x,t)\in \OOO \times I.
$$
Let $v \in H^1(0,T;H^1(\OOO)) \cap L^2(0,T;H^2(\OOO))$ satisfy
$$
\ppp_tv - \Delta v - \sum_{j=1}^n b_j(x)\ppp_{x_j}v - c(x)v
= F(x,t), \quad x \in\OOO, \, 0<t<T.               \eqno{(2.6)}
$$

Then
\\
{\bf Lemma 2 (Carleman estimate for parabolic equation).}\\
{\it 
Let $\la>0$ be chosen sufficiently large and let $\beta > 0$ be
arbitrarily fixed.  Then there exist constants $s_0>0$ and 
$C>0$ such that 
\begin{align*}
& \int_{Q_I} \left\{ \frac{1}{s}\left(\vert \ppp_tv\vert^2
+ \sum_{i,j=1}^n \vert \ppp^2_{x_jx_i} v\vert^2\right)
+ s\vert \nabla v\vert^2 + s^3\vert v\vert^2\right\}\tilweight dxdt\\
\le& C\int_{Q_I} \vert F\vert^2 \tilweight dxdt
+ Cs^3\int_{\ppp\OOO\times I} (\vert \nabla_{x,t}v\vert^2
+ \vert v\vert^2) \tilweight d\Sigma\\
+& Cs^3\int_{\OOO} (\vert \nabla v(x,t_0+\delta)\vert^2
+ \vert v(x,t_0+\delta)\vert^2 + \vert \nabla v(x,t_0-\delta)\vert^2
+ \vert v(x,t_0-\delta)\vert^2) e^{2s\www{\va}(x,t_0+\delta)}dx
\end{align*}
for all $s\ge s_0$.}
\\

This is a classical Carleman estimate and we can prove similarly for example
to Lemma 7.1 in \cite{BY2017} or Theorem 3.2 in \cite{Ya} by keeping
all the boundary integrals of $v(\cdot, t_0\pm \delta)$ and
$v$ on $\ppp\OOO \times I$ which are produced in the proof.

The Carleman estimate Lemma 1 needs 
extra data $u(\cdot,-T)$ and $u(\cdot,T)$ of the solution,
while Lemma 2 requires such data not only at the end points of the time 
interval but also on $\ppp\OOO \times I$.   
In applying them to inverse problems, we can control these terms by the 
weight $e^{2s\va}$ or $e^{2s\www{\va}}$ because the functions $\va$ and 
$\www{\va}$ take smaller values on such subboundaries.  
This is the essence of our argument without the cut-off.
\section{Proof of Theorem 1.}

{\bf First Step.}\\
By (1.6) and (1.10) and $R \in H^1(0,T; L^{\infty}(\OOO)) \subset 
C([0,T]; L^{\infty}(\OOO))$ by (1.8), setting $y = \ppp_tu$,
we have 
$$
\left\{ \begin{array}{rl}
& \ppp_t^2y - \Delta y - \sum_{j=1}^n b_j(x)\ppp_{x_j}y - c(x)y
= \ppp_tR(x,t)f(x) \in L^2(Q), \quad x \in\OOO, \, 0<t<T, \\
& y(\cdot,0) = 0, \quad \ppp_ty(\cdot,0) = R(\cdot,0)f \quad \mbox{in $\OOO$},\\
& y\vert_{\ppp\OOO\times (0,T)} = 0.
\end{array}\right.                         \eqno{(3.1)}
$$
For the application of the Carleman estimate, we extend $y$ to $t \in (-T,T)$
by the odd extension: $y(\cdot,-t) = -y(\cdot,t)$ for $0 < t < T$, and 
we make the odd extension of $(\ppp_tR)(\cdot,t)f$ to $(-T,0)$.
Then, by $y(\cdot,0) = 0$ in $\OOO$, we can directly verify that 
$y \in H^2(-T,T;L^2(\OOO)) \cap L^2(-T,T;H^2(\OOO) \cap H^1_0(\OOO))$ and
$$
\ppp_t^2y - \Delta y - \sum_{j=1}^n b_j(x)\ppp_{x_j}y - c(x)y
= \ppp_tR(x,t)f(x) \quad \mbox{in $Q_{(-T,T)}$}.      
%       \eqno{(3.2)}
$$

We set 
$$
d_0 := \min_{x\in\ooo{\OOO}} \vert x-x_0\vert, \quad 
d_1 := \max_{x\in\ooo{\OOO}} \vert x-x_0\vert.        \eqno{(3.2)}
$$
We define $\va$ in $Q_{(-T,T)}$ by (2.2) with $t_0=0$.
Since (1.9) means $T>\sqrt{d_1^2-d_0^2}$, 
we can choose $\beta \in (0,1)$ sufficiently close to $1$, such that 
$$
T > \frac{\sqrt{d_1^2-d_0^2}}{\sqrt{\beta}}.              \eqno{(3.3)}
$$
Therefore we can apply Lemma 1 to $y$ in $Q_{(-T,T)}$:
\begin{align*}
& \int_{Q_{(-T,T)}} (s\vert \nabla_{x,t}y\vert^2 + s^3\vert y\vert^2) 
\weight dxdt\\
\le& C\int_{Q_{(-T,T)}} \vert\ppp_tR\vert^2\vert f\vert^2 \weight dxdt 
+ Ce^{Cs}\int_{\Gamma\times (-T,T)} \vert \ppp_{\nu}y\vert^2 d\Sigma\\
+& Cs^3\int_{\OOO} (\vert \nabla_{x,t} y(x,T)\vert^2
+ \vert y(x,T)\vert^2 
+ \vert \nabla_{x,t} y(x,-T)\vert^2 + \vert y(x,-T)\vert^2)e^{2s\va(x,T)}dx
\end{align*}
for all $s>s_0$.  

We recall that $Q = \OOO\times (0,T)$.
Noting that $y(\cdot,-t) = -y(\cdot,t)$ for $-T<t<T$, we obtain
$$
\int_Q (s\vert \nabla_{x,t}y\vert^2 + s^3\vert y\vert^2) \weight dxdt
                                      \eqno{(3.4)}
$$
$$
\le CJ 
+ Ce^{Cs}\int_{\Gamma\times (0,T)} \vert \ppp_t\ppp_{\nu}u\vert^2 d\Sigma
+ Cs^3\int_{\OOO} (\vert \nabla_{x,t} y(x,T)\vert^2
+ \vert y(x,T)\vert^2) e^{2s\va(x,T)}dx
$$
for all $s\ge s_0$. Henceforth we set 
$$
J:= \int_{Q} \vert\ppp_tR\vert^2 \vert f\vert^2 \weight dxdt.
$$
\\
{\bf Second Step.}\\
We prove that there exist $s_0 > 0$ and $C>0$ such that 
$$
\int_{\Omega} |\ppp_t y(x,0)|^2e^{2s\va(x,0)} dx       \eqno{(3.5)}
$$
$$ 
\le CJ + C\int_{Q} (s|\nabla_{x,t}y|^2 + \vert y\vert^2)
e^{2s\va} dxdt + \int_\Omega |\nabla_{x,t} y(x,T)|^2e^{2s\va(x,T)} dx
$$
for all $s \ge 0$.
\\
{\bf Proof of (3.5).}\\
By direct calculations we can prove as follows. 
$$
\int_{\Omega} |\ppp_t y(x,0)|^2e^{2s\va(x,0)} dx 
= - \int_Q \ppp_t (|e^{s\va}\ppp_t y|^2) dxdt 
+ \int_{\Omega} |\ppp_t y(x,T)|^2e^{2s\va(x,T)} dx     \eqno{(3.6)}
$$
\begin{align*}
=& - \int_Q \left(2s(\partial_t \va)|\partial_t y|^2 
+ 2(\partial_t^2 y)\partial_t y\right) e^{2s\va} dxdt 
+ \int_{\Omega} |\ppp_t y(x,T)|^2e^{2s\va(x,T)} dx\\
= &- 2\int_Q \left\{s(\partial_t \va)|\partial_t y|^2 
+ \ppp_ty\left(\Delta y + \sum_{j=1}^n b_j\ppp_{x_j}y 
+ cy + (\ppp_tR)f\right) \right\} e^{2s\va} dxdt \\
+ & \int_{\Omega} |\ppp_t y(x,T)|^2e^{2s\va(x,T)} dx. 
\end{align*}
By $y = 0$ on $\ppp\OOO\times (0,T)$ and $y(\cdot,0) = 0$ in 
$\OOO$, integrating by parts, we estimate the following integral 
on the right-hand side in terms of (3.6): 
\begin{align*}
&-2 \int_Q \ppp_t y \Delta y e^{2s\va} dxdt 
= 2\int_Q \left(\nabla (\ppp_t y) \cdot\nabla y 
+ 2s(\partial_t y)\nabla\va\cdot\nabla y \right)e^{2s\va} dxdt\\
= & \int_\Omega |\nabla y(x,T)|^2e^{2s\va(x,T)} dx
- 2s\int_Q \partial_t\varphi\vert \nabla y\vert^2e^{2s\va} dxdt
+ 4\int_Q  s(\partial_t y)\nabla\va\cdot\nabla y e^{2s\va} dxdt.
\end{align*}
Since
$$
\vert (\ppp_ty)\nabla\va \cdot \nabla y\vert 
\le C(\vert \ppp_ty\vert^2 + \vert \nabla y\vert^2) \quad 
\mbox{in $Q$}
$$
and 
$$
\vert \ppp_ty\vert 
\left\vert \sum_{j=1}^n b_j\ppp_{x_j}y + cy + (\ppp_tR)f\right\vert 
\le C(\vert \nabla_{x,t} y\vert^2 + \vert y\vert^2 + 
\vert \ppp_tR\vert^2\vert f\vert^2) \quad \mbox{in $Q$},
$$ 
with (3.6) we can complete the proof of (3.5).
\\
{\bf Third Step.}\\
We will complete the proof of Theorem 1 by (3.4) and (3.5). 
The second equation in (3.1) implies
$$
\ppp_ty(\cdot,0) = R(\cdot,0) f \quad \mbox{in }\Omega.
$$
Therefore, by noting the assumption $\vert R(x,0)\vert \ne 0$ for
$x \in \ooo{\OOO}$ by (1.8), estimate (3.5) yields
\begin{align*}
& \int_{\Omega} |f(x)|^2e^{2s\va(x,0)} dx \\
\le & CJ + C\int_Q (s|\nabla_{x,t}y|^2 + \vert y\vert^2)e^{2s\va} dxdt 
+ \int_\Omega |\nabla_{x,t} y(x,T)|^2e^{2s\va(x,T)}dx.
\end{align*}
Applying (3.4) to the second term on the right-hand side to obtain
$$
\int_{\OOO} \vert f(x)\vert^2 e^{2s\va(x,0)} dx 
\le CJ + Ce^{Cs}\|\ppp_t\ppp_\nu u\|_{L^2(\Gamma\times(0,T))}^2
                                                 \eqno{(3.7)}
$$
$$
+ Cs^3\int_{\Omega} (|\nabla_{x,t} y(x,T)|^2 + |y(x,T)|^2)e^{2s\va(x,T)} dx
$$
for sufficiently large $s>0$. 

On the other hand, we have
$$
J = o(1)\int_\Omega |f(x)|^2 e^{2s\va(x,0)} dx \quad 
\mbox{as $s \to \infty$}.                \eqno{(3.8)}
$$
Indeed 
$$
e^{-2s(\va(x,0) - \va(x,t))}
= e^{-2se^{\lambda |x-x_0|^2}(1-e^{-\lambda \beta t^2})}
\le e^{-2s(1-e^{-\lambda \beta t^2})},
$$
by $e^{\lambda |x-x_0|^2} \ge 1$ for $x \in \ooo{\OOO}$, and so 
\begin{align*}
&J \le C\int_{\OOO} \vert f(x)\vert^2 e^{2s\va(x,0)}\left(
\int^T_0 \Vert\ppp_tR(\cdot,t)\Vert_{L^{\infty}(\OOO)}^2
e^{-2s(\va(x,0)-\va(x,t))}dt \right) dx\\
\le &C\int_{\OOO} \vert f(x)\vert^2 e^{2s\va(x,0)}\left(
\int^T_0 \Vert\ppp_tR(\cdot,t)\Vert_{L^{\infty}(\OOO)}^2
e^{-2s(1-e^{-\lambda \beta t^2})} dt \right) dx.
\end{align*}
Since $e^{-2s(1-e^{-\lambda \beta t^2})} \longrightarrow 0$ as 
$s\to \infty$ for fixed $0 < t \le T$ and
$\Vert\ppp_tR(\cdot,t)\Vert_{L^{\infty}(\OOO)}^2 \in L^1(0,T)$, 
we apply the Lebesgue convergence theorem, so that we can verify (3.8).

Therefore we absorb the first term on the right-hand side of (3.7) into the 
left-hand side:
$$
\int_\Omega |f(x)|^2 e^{2s\va(x,0)} dx 
                                               \eqno{(3.9)}
$$
$$
\le Cs^3\int_{\Omega} (|\nabla_{x,t} y(x,T)|^2 + |y(x,T)|^2)e^{2s\va(x,T)} dx 
+ Ce^{Cs}\|\ppp_t\ppp_\nu u\|_{L^2(\Gamma\times(0,T))}^2
$$
for sufficiently large $s$. 
Here we apply the classical a priori estimate 
(e.g., Lions and Magenes \cite{LM}) to (3.1), and we see
$$
\int_{\OOO} \vert \nabla_{x,t}y(x,T)\vert^2 dx 
\le C\Vert f\Vert^2_{L^2(\OOO)}.
$$
Moreover the Poincar\'e inequality yields
$$
\int_{\OOO} \vert y(x,T)\vert^2 dx 
\le C\int_{\OOO} \vert \nabla y(x,T)\vert^2 dx. 
$$ 
Hence
\begin{align*}
&\int_{\Omega} (|\nabla_{x,t} y(x,T)|^2 + |y(x,T)|^2)e^{2s\va(x,T)} dx \\
\le &Ce^{2se^{\lambda (d_1^2 - \beta T^2)}}
\int_\Omega (|\nabla_{x,t} y(x,T)|^2 + |y(x,T)|^2) dx
\le Ce^{2se^{\lambda (d_1^2 - \beta T^2)}}\|f\|_{L^2(\Omega)}^2. 
\end{align*}
On the other hand, we have
\begin{align*}
\int_\Omega |f(x)|^2 e^{2s\va(x,0)} dx 
= \int_\Omega e^{2se^{\lambda |x-x_0|^2}} |f(x)|^2 dx 
\ge e^{2se^{\lambda d_0^2}}\|f\|_{L^2(\Omega)}^2.
\end{align*}
Consequently (3.9) yields 
$$
\|f\|_{L^2(\Omega)}^2 \le Cs^3 e^{-c_0 s} \|f\|_{L^2(\Omega)}^2 
+ Ce^{Cs}\| \ppp_t\ppp_\nu u\|_{L^2(\Gamma\times(0,T))}^2.
$$
We set $c_0 = 2\left(e^{\lambda d_0^2} - e^{\lambda d_1^2 
- \lambda\beta T^2}\right)$.  The inequality (3.3) yields 
$c_0 > 0$.
Finally, by noting $\lim_{s\rightarrow \infty} s^3 e^{-c_0s} = 0$, 
we can absorb the first term on the right-hand side by taking sufficiently 
large $s$. This proves Theorem 1. $\blacksquare$
\section{Proof of Theorem 2}

{\bf First Step.}\\
We recall that $u \in H^2(0,T;H^1(\OOO)) \cap H^1(0,T;H^2(\OOO))$ satisfies 
(1.11).
Setting $z = \ppp_tu$ , we have
$$
\ppp_tz - \Delta z - \sum_{j=1}^n b_j(x)\ppp_{x_j}z - c(x)z
= \ppp_tR(x,t)f(x), \quad (x,t) \in Q_I           \eqno{(4.1)}
$$
and
$$
z(x,t_0) = \Delta u(x,t_0) + \sum_{j=1}^n b_j\ppp_{x_j}u(x,t_0) + cu(x,t_0) + R(x,t_0)f(x),
\quad x \in \OOO.                             \eqno{(4.2)}
$$
We apply Lemma 2 to $z$, and we obtain
$$
\int_{Q_I} \left( \frac{1}{s}\vert \ppp_tz\vert^2
+ s^3\vert z\vert^2\right) \tilweight dxdt              \eqno{(4.3)}
$$
$$
\le C\www{J} + Cs^3\int_{\ppp\OOO \times I} (\vert \nabla_{x,t}z\vert^2
+ \vert z\vert^2) \tilweight d\Sigma
$$
$$
+ Cs^3\int_{\OOO} (\vert \nabla z(x,t_0+\delta)\vert^2 
+ \vert z(x,t_0+\delta)\vert^2  
+ \vert \nabla z(x,t_0-\delta)\vert^2 
+ \vert z(x,t_0-\delta)\vert^2) e^{2s\www{\va}(x,t_0+\delta)} dx
$$
$$
=: \mathcal J_1+\mathcal J_2+\mathcal J_3.
$$
Here we set 
$$
\www{J} := \int_{Q_I} \vert \ppp_tR\vert^2\vert f(x)\vert^2
\tilweight dxdt.
$$
By (1.13) and the trace theorem, dividing the integral in (4.3) over
$\ppp\OOO \times I$ into $\Gamma\times I$ and 
$(\ppp\OOO\setminus \Gamma) \times I$, we can estimate
$$
\vert\mathcal J_2\vert+\vert\mathcal J_3\vert
                                            \eqno{(4.4)}
$$
\begin{align*}
\le & Ce^{Cs}\int_{\Gamma\times I} (\vert \nabla_{x,t}\ppp_tu\vert^2
+ \vert \ppp_tu\vert^2) d\Sigma
+ Cs^3M^2\exp\left( 2s\max_{x\in \ooo{\ppp\OOO\setminus \Gamma},
t \in \ooo{I}} \www{\va}(x,t)\right)\\
+ & Cs^3M^2 \exp\left( 2s\max_{x\in \ooo{\OOO}}
\www{\va}(x,t_0-\delta)\right).
\end{align*}
Since $\min_{x\in \ooo{\OOO_0}} d(x) > 0$ by (2.5), 
for $\delta > 0$, we can choose sufficiently large 
$\beta > 0$ such that $\max_{x\in \ooo{\OOO}} d(x) - \beta\delta^2
< 0$, and  
$$
\sigma_1:= \max\{ 
\max_{x\in\ooo{\ppp\OOO\setminus \Gamma}, t\in \ooo{I}} \www{\va}(x,t)
, \, \max_{x\in\ooo{\OOO}} \www{\va}(x,t_0-\delta) \}
< \sigma_0 := \min_{x\in \ooo{\OOO_0}} \www{\va}(x,t_0).
                                          \eqno{(4.5)}
$$
Indeed (4.5) is equivalent to 
$$
\max\{ \max_{x\in\ooo{\ppp\OOO\setminus \Gamma}} d(x), 
\, \max_{x\in\ooo{\OOO}} d(x) - \beta \delta^2\} 
< \min_{x\in \ooo{\OOO_0}} d(x).
$$
Since $d(x) = 0$ for $x \in \ppp\OOO\setminus \Gamma$ by 
(2.5), we can verify (4.5).
\\
         
Hence (4.4) yields
$$
\vert\mathcal J_2\vert+\vert\mathcal J_3\vert
\le Ce^{Cs}D^2 + Cs^3M^2e^{2s\sigma_1},
$$
where we set 
$$
D = \Vert \nabla_{x,t}\ppp_tu\Vert_{L^2(\Gamma\times I)}
+ \Vert \ppp_tu\Vert_{L^2(\Gamma\times I)}.
$$
Consequently (4.3) implies 
$$
\int_{Q_I} \left( \frac{1}{s}\vert \ppp_tz\vert^2
+ s^3\vert z\vert^2\right) \tilweight dxdt    
\le C\www{J} + Cs^3M^2e^{2s\sigma_1} + Ce^{Cs}D^2
                                                  \eqno{(4.6)}
$$
for all $s\ge s_0$.
\\
{\bf Second Step.}\\
We have
\begin{align*}
& \int_{\OOO} \vert z(x,t_0)\vert^2 e^{2s\www{\va}(x,t_0)} dx\\
= &\int^{t_0}_{t_0-\delta} \left( \ppp_t\int_{\OOO}
\vert z(x,t)\vert^2e^{2s\www{\va}(x,t)} dx \right) dt 
+ \int_{\OOO} \vert z(x,t_0-\delta)\vert^2 e^{2s\www{\va}(x,t_0-\delta)}dx\\
=& \int^{t_0}_{t_0-\delta} \int_{\OOO}
(2z\ppp_tz + 2s(\ppp_t\www{\va})\vert z\vert^2) e^{2s\www{\va}(x,t)} dxdt
+ \int_{\OOO} \vert z(x,t_0-\delta)\vert^2 e^{2s\www{\va}(x,t_0-\delta)}dx.
\end{align*}
Therefore, applying (4.5) to the second term on the right-hand side,
we obtain  
$$
\int_{\OOO} \vert z(x,t_0)\vert^2 e^{2s\www{\va}(x,t_0)} dx
\le C\int_{Q_I} (\vert z\vert \vert \ppp_tz\vert  + s\vert z\vert^2) 
e^{2s\www{\va}(x,t)} dxdt
+ CM^2 e^{2s\sigma_1}.                        \eqno{(4.7)}
$$
For the final term, we used (1.13).
Since 
$$
\vert z\vert\vert \ppp_tz\vert = s\vert z\vert \frac{1}{s}\vert \ppp_tz\vert
\le \frac{1}{2}\left(s^2\vert z\vert^2+ \frac{1}{s^2}\vert \ppp_tz\vert^2
\right),
$$
applying (4.6) and (4.7), we reach 
$$
\int_{\OOO} \vert z(x,t_0)\vert^2 e^{2s\www{\va}(x,t_0)}dx
\le \frac{C}{s}\www{J} + Cs^2M^2e^{2s\sigma_1} + Ce^{Cs}D^2
                                                      \eqno{(4.8)}
$$
for all $s\ge s_0$.  By (4.2) and the second condition in $(1.12)$, 
we estimate
\begin{align*}
& \int_{\OOO} \vert z(x,t_0)\vert^2 e^{2s\www{\va}(x,t_0)}dx\\
\ge& \int_{\OOO} \vert R(x,t_0)f(x)\vert^2 e^{2s\www{\va}(x,t_0)}dx
- C\int_{\OOO} \left\vert 
\Delta u(x,t_0) + \sum_{j=1}^n b_j\ppp_{x_j}u(x,t_0) 
+ cu(x,t_0)\right\vert^2 e^{2s\www{\va}(x,t_0)} dx \\
\ge& r_0^2\int_{\OOO} \vert f(x)\vert^2 e^{2s\www{\va}(x,t_0)}dx
- Ce^{Cs}\Vert u(\cdot,t_0)\Vert^2_{H^2(\OOO)}.
\end{align*}
Hence (4.8) yields
$$
\int_{\OOO} \vert f(x)\vert^2 e^{2s\www{\va}(x,t_0)}dx
\le C\www{J} + Cs^2M^2e^{2s\sigma_1} + Ce^{Cs}\www{D}^2, 
                                                \eqno{(4.9)}
$$
where we set $\www{D} = D + \Vert u(\cdot,t_0)\Vert_{H^2(\OOO)}$.

Since 
$$
\www{J} 
\le \int_{\OOO} \vert f(x)\vert^2 e^{2s\www{\va}(x,t_0)}
\left( \int^{t_0+\delta}_{t_0-\delta} \Vert \ppp_tR(\cdot,t)\Vert
_{L^{\infty}(\OOO)}^2 e^{-2s(\www{\va}(x,t_0) - \www{\va}(x,t))} dt
\right) dx,
$$
similarly to (3.8), we can verify 
$$ 
\www{J}  
= o(1) \int_{\OOO} \vert f(x)\vert^2 e^{2s\www{\va}(x,t_0)} dx 
\quad\mbox{as}\,\,\,s \to \infty.
$$
Therefore (4.9) implies 
$$
(1-o(1)) \int_{\OOO} \vert f(x)\vert^2 e^{2s\www{\va}(x,t_0)} dx
\le Cs^2M^2e^{2s\sigma_1} + Ce^{Cs}\www{D}^2\quad \forall s\ge s_0.
$$
Shrinking the integral domain $\OOO$ to $\OOO_0$
and using $\sigma_0 = \min_{x\in \ooo{\OOO_0}} \www{\va}(x,t_0)$, we see 
$$
\int_{\OOO_0} \vert f(x)\vert^2 dx e^{2s\sigma_0} 
\le Cs^2M^2e^{2s\sigma_1} + Ce^{Cs}\www{D}^2,
$$
that is,
$$
\Vert f\Vert^2_{L^2(\OOO_0)} \le Cs^2M^2e^{-2s\mu} + Ce^{Cs}\www{D}^2 \quad \forall s\ge s_0,
$$ 
where we have 
$$
\mu:= \sigma_0 - \sigma_1 > 0
$$
by (4.5).
Since $\sup_{s>0} s^2e^{-s\mu} < \infty$, replacing $C>0$ by $Ce^{Cs_0}$
and changing $s$ into $s+s_0$ with $s\ge 0$, we obtain
$$
\Vert f\Vert^2_{L^2(\OOO)} \le CM^2e^{-s\mu} + Ce^{Cs}\www{D}^2  
\quad \forall s\ge 0.                            \eqno{(4.10)}
$$
We minimize the right-hand side by choosing an appropriate value of 
parameter $s \ge 0$.
\\
{\bf Case 1: $M^2 > \www{D}^2$.} Then we can solve 
$$
M^2e^{-s\mu} = e^{Cs}\www{D}^2, \quad \mbox{that is,}\quad
s = \frac{2}{C+\mu}\log \frac{M}{\www{D}} > 0,
$$
so that 
$$
\Vert f\Vert^2_{L^2(\OOO_0)} \le CM^{2(1-\theta)}\www{D}^{2\theta},
$$
where $\theta = \frac{\mu}{C+\mu} \in (0,1)$.
\\
{\bf Case 2: $M^2 \le \www{D}^2$.} 
Then $\Vert f\Vert^2_{L^2(\OOO_0)} \le C(1 + e^{Cs})\www{D}^2$.
By the trace theorem and the Sobolev embedding, we readily see that 
$\www{D} \le CM$, and $\www{D} = \www{D}^{\theta}\www{D}^{1-\theta}
\le (CM)^{1-\theta}\www{D}^{\theta}$.  

Therefore, in both Cases 1 and 2,  we can obtain
$$
\Vert f\Vert^2_{L^2(\OOO_0)} \le C(M)\www{D}^{2\theta}.
$$
Thus the proof of Theorem 2 is completed. $\blacksquare$
\section{Concluding Remarks}

\mbox{\bf 5-1.} 
The method by Carleman estimates is widely applicable to other problems,
and as such a problem, we establish 
\\
{\bf Proposition 1 (observability inequality):}\\
For arbitrarily fixed $x_0 \not\in \ooo{\OOO}$, we assume (1.7) and 
$$
T > 2\sqrt{\max_{x\in \ooo{\OOO}} \vert x-x_0\vert^2
- \min_{x\in \ooo{\OOO}} \vert x-x_0\vert^2}.    \eqno{(5.1)}
$$
Then there exists a constant $C>0$ such that 
$$
\Vert u(\cdot,0)\Vert_{H^1_0(\OOO)} + \Vert \ppp_tu(\cdot,0)\Vert
_{L^2(\OOO)} \le C\Vert \ppp_{\nu}u\Vert_{L^2(\Gamma \times (0,T))}
                                          \eqno{(5.2)}
$$
for each $u$ satisfying
$$
\left\{ \begin{array}{rl}
& \ppp_t^2u = \Delta u + \sum_{j=1}^n b_j(x)\ppp_{x_j}u + c(x)u \quad
\mbox{in $Q$}, \\
& u(\cdot,0) \in H^1_0(\OOO), \quad \ppp_tu(\cdot,0) \in L^2(\OOO),\\
& u\vert_{\ppp\OOO\times (0,T)} = 0,        
\end{array}\right.                      \eqno{(5.3)}
$$
where $b_j, c \in L^{\infty}(\OOO)$, $j=1,..., n$.
\\

Inequality (5.2) is called an observability inequality, and  
there are many related works in the control theory (e.g., \cite{Ko}).  
The proof by Carleman estimates is found 
for example, in Chapter 4 in \cite{BY2017}, \cite{FLZ}, pp.58-65 in \cite{KlT}.
Our proposed argument in this article can simplify the 
existing proofs, as one sees below.  

%Before the proof, we make comparisons with some of the existing results.
By the finiteness of the propagation speed for the hyperbolic equation,
the observation time $T$ cannot be arbitrary for estimate (5.2).
The right-hand side of (5.1) gives a critical value of $T$, which 
can be described only by a choice of $x_0$ and $\OOO$.
Other papers give different critical values and we can compare 
for example, formula (4.28) (p.96) in \cite{FLZ} and formula (14) (p.36) 
in \cite{Ko}, which are worse than ours (5.1) for the case
of $\ppp_t^2 - \Delta$ as the principal term of the hyperbolic 
equation.
In (5.1), we do not consider the case $x_0\in \ooo{\OOO}$, but we 
can similarly discuss also for the case of $x_0 \in \ooo{\OOO}$.
Here we omit the discussions for showing the essence of our method.

{\bf Proof.}\\
We recall $Q = \OOO \times (0,T)$ and 
(3.2): $d_0 = \min_{x\in \ooo{\OOO}} \vert x-x_0\vert$ and
$d_1 = \max_{x\in \ooo{\OOO}} \vert x-x_0\vert$, and we set 
$$
\kappa_1 = \exp\left( \lambda \left(d_1^2 - \frac{T^2}{4}\beta
\right)\right), \quad
\kappa_0 = \exp(\lambda d_0^2).
$$
We replace the time interval $(-T,T)$ by $(0,T)$ and we apply 
Lemma 1 in $Q:= \OOO \times (0,T)$.
We choose $t_0 = \frac{T}{2}$, and (5.1) allows us to take 
$0 < \beta < 1$ in (2.2) such that 
$$
T > 2\frac{\sqrt{d_1^2 - d_0^2}}{\sqrt{\beta}}.  \eqno{(5.4)}
$$
Then $d_1^2 - \frac{T}{4}\beta^2 < d_0^2$, that is,
$\kappa_0 > \kappa_1$.

Now we employ Lemma 1 to (5.3): 
\begin{align*}
& \int_Q s\vert \nabla_{x,t}u\vert^2 \weight dxdt
\le Ce^{Cs}\Vert \ppp_{\nu}u\Vert^2_{L^2(\Gamma\times (0,T))}\\
+ &C\int_{\OOO} (s\vert \nabla_{x,t}u(x,0)\vert^2
+ s^3\vert u(x,0)\vert^2 
+ s\vert \nabla_{x,t}u(x,T)\vert^2 + s^3\vert u(x,T)\vert^2)
e^{2s\va(x,0)} dx
\end{align*}
for all large positive $s.$
We set $E(t) = \int_{\OOO} \vert \nabla_{x,t}u(x,t)\vert^2 dx$ for
$0\le t \le T$.
Then the classical energy estimate (e.g., \cite{LM}) and 
the Poincar\'e inequality yield
$$
\int_{\OOO} (s\vert \nabla_{x,t}u(x,0)\vert^2
+ s^3\vert u(x,0)\vert^2 
+ s\vert \nabla_{x,t}u(x,T)\vert^2 + s^3\vert u(x,T)\vert^2)
e^{2s\va(x,0)} dx \le Cs^3E(0)e^{2s\kappa_1}.
$$
Hence
$$
\int_Q s\vert \nabla_{x,t}u\vert^2 \weight dxdt 
\le Ce^{Cs}\Vert \ppp_{\nu}u\Vert^2_{L^2(\Gamma\times (0,T))}
+ Cs^3e^{2s\kappa_1}E(0).                   \eqno{(5.5)}
$$
By (5.4) we further find small $\delta > 0$ such that 
$T > 2\frac{\sqrt{d_1^2 - d_0^2 + \beta \delta^2}}{\sqrt{\beta}}$.
Then we can directly verify
$$
\kappa_2:= e^{\lambda(d_0^2 - \beta \delta^2)} > \kappa_1.
                                       \eqno{(5.6)}
$$
Hence, since $\va \ge \kappa_2$ on $\ooo{\OOO} \times \left[ 
\frac{T}{2} - \delta, \frac{T}{2} + \delta\right]$, we obtain
$$
\int_Q s\vert \nabla_{x,t}u\vert^2 \weight dxdt
\ge \int_{\frac{T}{2}-\delta}^{\frac{T}{2}+\delta}\int_{\OOO}
s\vert \nabla_{x,t}u\vert^2 \weight dxdt
\ge se^{2s\kappa_2}\int_{\frac{T}{2}-\delta}^{\frac{T}{2}+\delta}
E(t) dt.
$$
Again with the classical energy estimate, this yields 
$$
\int_Q s\vert \nabla_{x,t}u\vert^2 \weight dxdt
\ge 2Cse^{2s\kappa_2}\delta E(0).
$$
Therefore (5.5) yields 
$$
2Cse^{2s\kappa_2}\delta E(0) \le Ce^{Cs}\Vert \ppp_{\nu}u\Vert^2
_{L^2(\Gamma\times (0,T))}
+ Cs^3e^{2\kappa_1}E(0),
$$
that is,
$$
2Cse^{2s\kappa_2}\delta \left( 1 - \frac{C_1}{\delta}s^2
e^{-2s(\kappa_2-\kappa_1)}\right)E(0)
\le Ce^{Cs}\Vert \ppp_{\nu}u\Vert^2_{L^2(\Gamma\times (0,T))}.
$$
By (5.6), choosing $s>0$ large, we complete the proof of the 
observability inequality.
\\

{\bf 5.2.}
Our method is applicable to a Cauchy problem for a parabolic equation.
\\
{\bf Proposition 2.}
\\
Let $\Gamma \subset \ppp\OOO$ and $\OOO_0 \subset \OOO$ be given as in 
Theorem 2.  We assume that $u \in H^1(0,T;L^2(\OOO)) \cap
L^2(0,T;H^2(\OOO))$ satisfy
$$
\ppp_tu = \Delta u + \sum_{j=1}^n b_j(x)\ppp_{x_j}u + c(x)u \quad
\mbox{in $Q$},                   \eqno{(5.7)}
$$
with $b_j, c \in L^{\infty}(\OOO)$, $j=1,..., n$, and 
$$
\Vert u\Vert_{H^1(0,T;L^2(\OOO))} + \Vert u\Vert
_{L^2(0,T;H^2(\OOO))} \le M                \eqno{(5.8)}
$$
with some constant $M>0$.  Let $\ep \in (0,T)$ be arbitrarily given.
Then there exist constants $C>0$ and $\theta \in (0,1)$ such that 
$$
\Vert u\Vert_{H^1(\ep,T-\ep;L^2(\OOO_0))} + \Vert u\Vert
_{L^2(\ep,T-\ep;H^2(\OOO_0))}
\le C(\Vert \nabla_{x,t}u\Vert_{L^2(\Gamma \times (0,T))}
+ \Vert u\Vert_{L^2(\Gamma \times (0,T))})^{\theta}.
$$
\\

This is a conditional stability estimate for the Cauchy problem for 
a parabolic equation (5.7) and see e.g., Theorem 5.1 in 
\cite{Ya}.  Our proof is much simpler.

{\bf Proof.}\\
For given $t_0>0$ and $\www{\delta} > 0$ satisfying $0<t_0-\www{\delta}
<t_0 + \wdel < T$.  We apply the Carleman estimate Lemma 2 in 
$Q_{(t_0-\wdel,t_0+\wdel)} := \OOO \times (t_0-\wdel, t_0+\wdel)$ to obtain
$$
\int_{Q_{(t_0-\wdel,t_0+\wdel)}} \left\{
\frac{1}{s}\left( \vert \ppp_tu\vert^2 + \sum_{i,j=1}^n
\vert \ppp^2_{x_ix_j}u\vert^2\right) + s\vert \nabla u\vert^2
+ s^3\vert u\vert^2\right\} e^{2s\www{\va}} dxdt 
                                           \eqno{(5.9)}
$$
\begin{align*}
\le &Cs^3\int_{\ppp\OOO \times (t_0-\wdel,t_0+\wdel)}
(\vert \nabla_{x,t}u\vert^2 + \vert u\vert^2) e^{2s\www{\va}} d\Sigma\\
+ &Cs^3\int_{\OOO} (\vert \nabla u(x, t_0-\wdel)\vert^2
+ \vert u(x,t_0-\wdel)\vert^2 + \vert \nabla u(x,t_0+\wdel)\vert^2
+ \vert u(x,t_0+\wdel)\vert^2) e^{2s\www{\va}(x,t_0+\wdel)} dx
\end{align*}
for all $s \ge s_0$.
Here we note that the constants $C>0$ and $s_0>0$ are independent of 
$t_0$ because the Carleman estimate is invariant by the translation in time
provided that the translated time interval is in $(0,T)$.

We set 
$$
\www{d_0} = \min_{x\in \ooo{\OOO_0}} d(x), \quad
\www{d_1} = \max_{x\in \ooo{\OOO}} d(x),
$$
and choose $\weps>0$ such that $0<\weps<\wdel$.  Then we have
$$
\max_{x\in \ooo{\ppp\OOO\setminus \Gamma},
t_0-\wdel\le t \le t_0+\wdel} \www{\va}(x,t) \le 1,
\quad \max_{x\in \ooo{\OOO}} \www{\va}(x,t_0-\wdel)
= \max_{x\in \ooo{\OOO}} \www{\va}(x,t_0+\wdel)
= e^{\lambda(\www{d_1}^2-\beta \wdel^2)}
$$
and
$$
\min_{x\in \ooo{\OOO_0}, t_0-\weps\le t\le t_0+\weps} \www{\va}(x,t)
\ge e^{\lambda(\www{d_0}^2-\beta \weps^2)}.
$$
Therefore, shrinking the integral domain $Q_{(t_0-\wdel,t_0+\wdel)}$
to $\OOO_0 \times (t_0-\weps,t_0+\weps)$ 
in the left-hand side of (5.9), by (5.8) we obtain
$$
\frac{1}{s}\exp(2se^{\lambda(\www{d_0}^2-\beta \weps^2)})
(\Vert u\Vert^2_{H^1(t_0-\weps,t_0+\weps;L^2(\OOO_0))}
+ \Vert u\Vert^2_{L^2(t_0-\weps,t_0+\weps;H^2(\OOO_0))})    \eqno{(5.10)}
$$
\begin{align*}
\le &Cs^3\left(\int_{\Gamma \times (t_0-\wdel,t_0+\wdel)}
(\vert \nabla_{x,t}u\vert^2 + \vert u\vert^2) e^{2s\www{\va}} d\Sigma
+ \int_{(\ppp\OOO \setminus \Gamma) \times (t_0-\wdel,t_0+\wdel)}
(\vert \nabla_{x,t}u\vert^2 + \vert u\vert^2) e^{2s\www{\va}} d\Sigma
\right)                             \\
+& Cs^3\int_{\OOO} (\vert \nabla u(x, t_0-\wdel)\vert^2
+ \vert u(x,t_0-\wdel)\vert^2 + \vert \nabla u(x,t_0+\wdel)\vert^2
+ \vert u(x,,t_0+\wdel)\vert^2) e^{2s\www{\va}(x,t_0+\wdel)} dx\\
\le& Cs^3e^{Cs}D^2 + Cs^3e^{2s}M^2
+ Cs^3M^2\exp(2se^{\lambda(\www{d_1}^2-\beta \wdel^2)})
\end{align*}
for all $s\ge s_0$  Here we set $D = \Vert \nabla_{x,t}u\Vert
_{L^2(\Gamma \times (0,T))} + \Vert u\Vert_{L^2(\Gamma \times (0,T))}$.

Now for given $\weps > 0$, we choose $\beta > 0$ and $\wdel > 0$.  
For $t_0 \in (\wdel, T-\wdel)$, we note $(t_0-\wdel, t_0+\wdel)
\subset (0,T)$. 
First choose large $N>1$ such that 
$$
N-1 > \frac{\www{d_1}^2 - \www{d_0}^2}{\www{d_0}^2},
$$
and set $\wdel = N\weps$.  Then, noting that $N^2-1 > N-1$, we can prove  
$$
\frac{\www{d_1}^2 - \www{d_0}^2}{\wdel^2 - \weps^2} 
< \frac{\www{d_0}^2}{\weps^2}.
$$
Therefore we can choose $\beta > 0$ such that 
$$
\frac{\www{d_1}^2 - \www{d_0}^2}{\wdel^2 - \weps^2} < \beta 
< \frac{\www{d_0}^2}{\weps^2}.
$$
With these chosen $\beta > 0$ and $\wdel > 0$, we can directly verify
$$
\mu_1:= e^{\lambda(\www{d_0}^2-\beta \weps^2)}
> \mu_2:= \max \{1, \, e^{\lambda(\www{d_1}^2-\beta \wdel^2)}\}.
$$
Hence (5.10) yields 
$$
\Vert u\Vert_{H^1(t_0-\weps,t_0+\weps;L^2(\OOO_0))}
+ \Vert u\Vert_{L^2(t_0-\weps,t_0+\weps;H^2(\OOO_0))}
\le Cs^4M^2e^{-2s\mu_0} + Cs^4e^{Cs}D^2
$$
for all $s\ge s_0$.  Here we note 
$$
\mu_0:= \mu_1 - \mu_2 > 0.
$$
Hence, arguing similarly to after (4.10), we obtain
$$
\Vert u\Vert_{H^1(t_0-\weps,t_0+\weps;L^2(\OOO_0))}
+ \Vert u\Vert_{L^2(t_0-\weps,t_0+\weps;H^2(\OOO_0))}
\le C(M)D^{\theta},
$$
where the constants $C(M)$ and $\theta \in (0,1)$ are dependent on 
$M, \weps, \wdel > 0$, but independent of $t_0$.
Varying $t_0$ over $(\wdel, T-\wdel)$, we have 
$$
\Vert u\Vert_{H^1(\wdel-\weps,T-\wdel+\weps;L^2(\OOO_0))}
+ \Vert u\Vert_{L^2(\wdel-\weps,T-\wdel+\weps;H^2(\OOO_0))}
\le C(M)D^{\theta}.
$$
For given $\ep>0$ in the statement of the proposition,
we choose $\weps = \frac{\ep}{N-1}$, so that $\wdel - \weps = (N-1)\weps
= \ep$ and $T-\wdel+\ep = T-\ep$, we can complete the proof of Proposition 2.  
\\

{\bf 5-3.} 
Our argument proposed in this article works for similar inverse problems
for various types of partial differential equations such as plate 
equations, Schr\"odinger equation, integro-differential equations,
Lam\'e equations, equations for fluid dynamics.

%%%%%%%%%%%%%%%%%%%%%%%%%%%%%%%%%%%%%%%%%%%%%%%%%
%%%%%%%%%%%%%%%%%%%%%%%%%%%%%%%%%%%%%%%%%%%%%%%%%

%%%%%%%%%%%%%%%%%%%%%%%%%%%%%%%%%%%%%%%%%%%%%%%%%
%%%%%%%%%%%%%%%%%%%%%%%%%%%%%%%%%%%%%%%%%%%%%%%%%
\section*{Acknowledgments}
The first author thanks the Leading Graduate Course for Frontiers of 
Mathematical Sciences and Physics (FMSP, The University of Tokyo) 
and was supported by Grant-in-Aid for Scientific Research (S) 
15H05740 and Grant-in-Aid for Research Activity Start-up 19K23400 of Japan 
Society for the Promotion of Science. 
The second author was partially supported by NSF grant DMS 1312900 and
Grant-in-Aid for Scientific Research (S) 
15H05740 of Japan Society for the Promotion of Science. 
The third author was supported by Grant-in-Aid for Scientific Research (S) 
15H05740 of Japan Society for the Promotion of Science and
by The National Natural Science Foundation of China 
(no. 11771270, 91730303).
This work was supported by A3 Foresight Program ``Modeling and Computation of 
Applied Inverse Problems" of Japan Society for the Promotion of Science and 
prepared with the support of the "RUDN University Program 5-100".

%%%%%%%%%%%%%%%%%%%%%%%%%%%%%%%%%%%%%%%%%%%%%%%%%
%%%%%%%%%%%%%%%%%%%%%%%%%%%%%%%%%%%%%%%%%%%%%%%%%

\end{document}